\documentclass[11pt]{article}

\usepackage{amsmath,amsfonts,amssymb}

\newcommand{\be}{\begin}
\newcommand{\Ref}[1]{(\ref{#1})}

\newcommand\ca{\mathcal A}
\newcommand{\cb}{\mathcal B}

\newcommand\cD{\mathcal D}

\newcommand{\cg}{\mathcal G}

\newcommand\cO{\mathcal O}

\newcommand\fl{\mathcal R}

\newcommand\ct{\mathcal T}

\newcommand{\al}{\alpha}
\newcommand{\del}{\delta}

\newcommand{\eps}{\epsilon}
\newcommand\ga{\gamma}
\newcommand\Ga{\Gamma}
\newcommand\ka{\kappa}
\newcommand\lla{\lambda}
\newcommand\La{\Lambda}
\newcommand\om{\omega}
\newcommand\Om{\Omega}
\newcommand{\si}{\sigma}
\newcommand\Si{\Sigma}

\newcommand{\co}{\mathbb C}

\newcommand{\R}{\mathbb R}

\newcommand{\z}{\mathbb Z}

\newcommand\CP[1]{{\mathbb{CP}^{#1}}}

\newcommand{\uset}{\underset}
\newcommand{\oset}{\overset}
\newcommand{\uline}{\underline}

\newcommand{\la}{\langle}
\newcommand{\ra}{\rangle}
\newcommand{\st}{\,|\,}
\newcommand{\ti}{\tilde}
\newcommand{\setmin}{\backslash}
\newcommand{\prtl}{\partial}
\renewcommand\square{\kern20pt{\vbox{\hrule height.4pt
        \hbox{\vrule width.4pt height 6pt\kern6pt
                \vrule width.4pt}
        \hrule height.4pt}}}

\newcommand\rank{\text{rank}}

\newcommand\inv{^{-1}}
\newcommand\hol{\text{Hol}}
\newcommand\aut{\text{Aut}}
\newcommand\waut{\widetilde{\aut}(E_z)}

\newcommand\rf{{\R}^4}
\newcommand\proof{{\em Proof.}\ }
\newcommand\loc{{\text{loc}}}

\newcommand{\ry}{\R\times Y}
\newcommand{\rpy}{\R_+\times Y}

\newcommand\U[1]{\text{U}(#1)}
\newcommand\SU[1]{\text{SU}(#1)}

\newcommand\SO[1]{\text{SO}(#1)}
\newcommand\so[1]{\text{so}(#1)}
\newcommand\Sp[1]{\text{Sp}(#1)}

\newcommand\CF{\text{CF}}
\newcommand\HF{\text{HF}}
\newcommand\ind{\text{ind}}

\newcommand\fd{\mathfrak d}

\newcommand\frakp{\mathfrak p}

\newcommand\lloc[2]{L^#1_{#2,\loc}}
\newcommand\llw[3]{{L^{#1,#2}_#3}}

\newcommand\lw[2]{{L^{#1,#2}}}

\newcommand\itr{\text{int}}

\newcommand\ck{\mathcal K}

\newcommand\aref{{A_{\text{ref}}}}

\renewcommand\O[1]{\text{O}(#1)}

\newcommand\tr{\text{tr}}
\newcommand\hola{\hol_\ga(A)}
\newcommand\holza{\hol_z(A)}

\newcommand\xg{\Xi_\gamma}
\def\zmt{{\z/2}}
\def\tc{\ti c_1}
\def\red{^{\text{red}}}
\def\tred{^{\text{tred}}}
\def\ared{^{\text{ared}}}
\def\mared{M\ared}

\def\tpl{\ti P_\ell}
\def\tpc{\ti P_c}

\def\ctl{\ct_\ell}

\def\hfl{\hat F_\ell}

\def\torsion{\text{torsion}}

\newcommand\hth{\hat\Theta}
\newcommand\hmk{\hat M_k}
\newcommand\dla{\uline\lla}

\begin{document}

\title{$4$--manifolds and intersection forms with\\local coefficients}

\author{Kim A.\ Fr\o yshov
\thanks{This work was supported by a ProDoc grant at the ETH,
  Z\"urich, and by QGM (Centre for Quantum Geometry of Moduli Spaces)
  funded by the Danish National Research Foundation.}}

\date{}

\maketitle

\be{abstract}
We extend Donaldson's diagonalization theorem to
intersection forms with certain local coefficients, under some
constraints. This provides new examples of non-smoothable topological
$4$--manifolds.
\end{abstract}

\thispagestyle{empty}

\bibliographystyle{plain}


\section{Introduction}

A celebrated early theorem of Donaldson \cite{D1,D2} says that if the
intersection form of a closed, oriented smooth $4$--manifold $V$ is negative
definite,
then it is {\em standard}, i.e.\ there is a basis for $H^2(V;\z)/\torsion$
with respect
to which the form is diagonal. The proof involved a careful study of a 
certain $\SU2$--instanton moduli space over $V$. Later, Fintushel and Stern
\cite{FS2}
found a simpler proof using $\SO3$--instanton moduli spaces in the
case when $H_1(V;\z)$ contains no $2$--torsion. (The assumption on the
torsion can be removed by using results from \cite{D2}, see \cite{Fr3}.)
In either variant of the proof
an essential point is the link between the intersection
form of $V$ and Abelian reducibles
in the moduli spaces, which are represented by connections with
stabilizer $\U1$. In $\SO3$--moduli spaces there is also 
a second type of reducible, namely the {\em twisted reducibles},
which are represented by connections with stabilizer $\z/2$
(among all automorphisms of the $\SO3$--bundle).
In this paper we will show that these are related
to the intersection forms of $V$ with certain local coefficients. We
use this to partially extend Donaldson's theorem to such forms.
We will now explain our result in more detail.

We generalize the setup somewhat and consider a compact, connected, oriented,
smooth $4$--manifold $X$ 
with boundary $Y$. Let $\ell\to X$ be any bundle of infinite cyclic groups.
Recall that the set of isomorphism classes of such bundles form an Abelian
group isomorphic to $H^1(X;\z/2)$. Let $H^*(X;\ell)$ be the singular
cohomology with $\ell$ as bundle of coefficients.
Since $\ell\otimes\ell=\z$,
the cup product defines a homomorphism
\be{equation}\label{eqn:cup-prod}
H^2(X;\ell)\otimes H^2(X,Y;\ell)\to H^4(X,Y;\z)=\z.
\end{equation}
Now suppose $Y$ is an integral homology sphere. Then
$H^2(X,Y;\ell)=H^2(X;\ell)$, and \Ref{eqn:cup-prod} induces a unimodular
quadratic form $Q_{X,\ell}$ on $H^2(X;\ell)/\text{torsion}$, which we refer to
as the intersection form of $X$ with coefficients in $\ell$. When $\ell$ is
trivial this is of course the usual intersection form of $X$.
The signature of $Q_{X,\ell}$ is independent of $\ell$.
As observed in \cite[p.\,587]{KM3}, the same holds for the quantity
\[b_0(X;\ell)-b_1(X;\ell)+b^+(X;\ell),\]
where $b_j(X;\ell):=\rank\, H^j(X;\ell)$ and
$b^+_2(X;\ell)$ denotes the dimension of a maximal positive subspace
for $Q_{X,\ell}$. For any non-trivial $\ell$ one therefore has
\be{equation}\label{eqn:del-indep}
-b_1(X;\ell)+b^+(X;\ell)=1-b_1(X)+b^+(X),
\end{equation}
where $b_j(X):=b_j(X;\z)$ and $b^+(X):=b^+(X;\z)$.

For any Abelian group $G$ let $\HF^*(Y;G)$ denote the instanton Floer 
cohomology group with coefficients in $G$, see \cite{F1,D5}. This is the 
cohomology of a cochain complex $\CF^*\otimes G$, where $\CF^q$ is the
free Abelian group generated by gauge equivalence classes of
irreducible (perturbed) flat $\SO3$--connections
over $Y$ of index $q\in\z/8$, and the differential
$d:\CF^q\to\CF^{q+1}$ counts instantons over the cylinder $\ry$
interpolating between two given irreducible flat connections.
Counting $\SO3$--instantons over $\ry$ with trivial flat limit at $+\infty$
yields a homomorphism $\del:\CF^4\to\z$ which satisfies $\del d=0$ 
(see \cite{Fr3}) and therefore induces a homomorphism
\[\del_0:\HF^4(Y;G)\to G.\]

Before stating the main result of this paper we need one more definition:
\be{align*}
\tau(X)&:=\dim_{\z/2}\left[\torsion(H_1(X;\z))\otimes\z/2\right]\\
       &=b_1(X;\z/2)-b_1(X),
\end{align*}
where $b_j(X;\z/2):=\dim_{\z/2}H_j(X;\z/2)$.

\be{thm}\label{thm:main-theorem}
Let $X$ be any compact, connected, oriented, smooth $4$--manifold
whose boundary $Y$ is
an integral homology sphere, and such that
\be{equation}\label{eqn:thmeqn}
\tau(X)+b^+(X)\le2.
\end{equation}
Let $\ell\to X$ be any non-trivial bundle of infinite cyclic groups. If
$Q_{X,\ell}$ is non-standard negative definite and
$H^2(X;\ell)$ contains no element of order $4$ then
\[\del_0:\HF^4(Y;\z/2)\to\z/2\]
is non-zero.
\end{thm}

\be{cor}\label{cor:Vclosed}
Let $V$ be any closed, connected, oriented, smooth $4$--manifold
such that
\[\tau(V)+b^+(V)\le2.\]
Let $\ell\to V$ be any non-trivial bundle of infinite cyclic groups such that
$Q_{V,\ell}$ is negative definite and
$H^2(V;\ell)$ contains no element of order $4$. Then $Q_{V,\ell}$ is standard.
\end{cor}

\proof This follows from the theorem by taking $X$ to be the complement
of an open $4$--ball in $V$, and recalling that $\HF^*(S^3;\z/2)=0$.\square

{\bf Remarks.} {\bf(i)}  Under
the hypotheses of the corollary, $V$ cannot be spin. For in that case
the usual intersection form $Q_V$ would be even with negative signature,
so $Q_V$ could not be definite by Donaldson's theorem. The condition
$b^+(V)\le2$ would then violate a theorem of Furuta \cite{Fur2}.

{\bf(ii)} If $b_1(X)=1$ and $\tau(X)=0$ then $H^2(X;\ell)$ does not even have
any element of order $2$, see Proposition~\ref{prop:order-two}.

{\bf(iii)} The author does not know whether the theorem holds without the
assumptions on $\tau(X)+b^+(X)$
and (in general) elements of order $4$ in $H^2(X;\ell)$, despite
attempts at finding counterexamples. 

{\bf(iv)} The statement of the theorem
holds when $\ell$ is trivial too, and without
the assumption $\tau(X)\le2$. 
However, we prefer to take that up in a separate paper.

{\bf(v)} One reason for the appearance
of the term $\tau+b^+$ in the theorem is that
this quantity is invariant under surgery on
any circle in the interior of $X$ which represents a non-zero
class in $H_1(X;\z/2)$, see Lemma~\ref{lemma:taubpl}.

\be{prop}
Let $V$ be any closed, oriented topological $4$--manifold whose
intersection form $Q_V$ is non-standard negative definite. Suppose $H_1(V;\z)$
contains no element of order $4$. Let either
\be{description}
\item[(i)]$W=\Si\times S^2$, where $\Si$ is any closed, oriented, connected surface of
genus at least $1$, or
\item[(ii)]$W=Y\times S^1$, where $Y$ is any closed, oriented $3$--manifold.
\end{description}
If $\tau(V)+\tau(W)+b^+(W)\le2$, then $V\#W$ does not admit any
smooth structure.
\end{prop}

Of course, if $W=\Si\times S^2$ then $\tau(W)=0$ and $b^+(W)=1$,
whereas if $W=Y\times S^1$ then $\tau(W)=\tau(Y)$ and $b^+(W)=b_1(Y)$.

\proof (i) We may assume that $V$ is connected and that $Q_V$ is
negative definite. Let $\ell\to W:=\Si\times S^2$
be any non-trivial bundle of infinite cyclic groups. The exact
sequence \Ref{eqn:les-double-cover} below yields
\[\torsion(H_1(W;\ell))=\z/2,\qquad H_2(W;\ell)=\z/2.\]
Let $\ell'\to V':=V\#W$ be the bundle which corresponds to the trivial bundle
over $V$ and to $\ell$ over $W$. Then the group
\[H_1(V';\ell')=H_1(V;\z)\oplus H_1(W;\ell)\]
contains no element of order $4$. By the universal coefficient theorem
(see \Ref{eqn:univ-coeff} below) the same holds for $H^2(V';\ell')$.
As for the intersection forms one has
\[Q_{V',\ell'}=Q_V,\]
so it follows from Corollary~\ref{cor:Vclosed} that $V'$ cannot admit any smooth
structure.

(ii) Let $\ell\to W:=Y\times S^1$ be the pull-back of the non-trivial
$\z$--bundle over $S^1$. Using the exact sequence
\Ref{eqn:les-double-cover} below one finds that $H_k(W;\ell)$ is a finite
group for all $k$, and that
\[H_1(W;\ell)=H_1(Y;\z)\,/\,2H_1(Y;\z)\approx(\z/2)^r\]
for some $r$. We can now argue as in (i).
\square \\

When combined with Freedman's classification of simply-connected, closed,
oriented
topological $4$--manifolds \cite{Freedman-Quinn}
this yields many examples of non-smoothable
indefinite $4$--manifolds, also with odd intersection form.
In the case of even intersection form
such examples can also be found using Rochlin's theorem
or Furuta's theorem.

Note that if $V$ is simply connected and negative definite, say, then
$V\#\CP2$ is smoothable, since by Freedman's theorem and
the classification of odd indefinite forms it is
homeomorphic to $\CP2\#(-n\CP2)$ for some $n$.

In a slightly different direction, Friedl--Hambleton--Melvin--Teichner
\cite{FHMT} have proved that a certain negative definite closed,
oriented topological $4$--manifold
$V$ with $\pi_1(V)=\z$ and $b_2(V)=4$ is not smoothable by applying
Donaldson's diagonalization theorem to the finite coverings of $V$.
(A survey of related material can be found in \cite{Hambleton1}.)

After some preliminaries in Section~\ref{sec:cohomology-local}
on (co)homology with local coefficients, Section~\ref{sec:conn-hol}
introduces what is probably the main novelty in the paper as far as
gauge theory is concerned: Given any $\SO3$--bundle $E\to Z$, where
$Z$ is a smooth, compact manifold, and any loop $\ga:S^1\to Z$, we define
a double covering $\Xi_\ga\to U_\ga$,
where $U_\ga$ is a certain open subset of the
orbit space $\cb(E)$ of all connections in $E$ (of a given Sobolev type).
The subset $U_\ga$ contains all irreducible connections
as well as some reducibles including all Abelian ones.
In Section~\ref{sec:tw-red} we classify non-flat
twisted reducible instantons over certain $4$--manifolds $W$ with a
tubular end. The local structure around these reducibles is described
in Section~\ref{sec:local-str-twr}, whereas Abelian reducibles
are discussed in
Section~\ref{sec:abelian-inst}. Section~\ref{sec:banach} proves
three lemmas on Banach manifolds.
Section~\ref{sec:proof-thm} contains the proof of the theorem. This begins by
reducing the problem to the case $b_1(X)=1+b^+(X)$ by doing surgery
on a suitable collection of disjoint circles in $X$.
We then study the moduli space $M_k$ of instantons with trivial limit
in a certain $\SO3$--bundle over $W:=X\cup_Y(\rpy)$. The irreducible part
$M^*_k$ is cut down to a $1$--manifold using sections of the real line
bundles corresponding to suitable double coverings $\Xi_\ga$.
The ends of this $1$--manifold are associated to twisted reducibles
in $M_k$ and factorizations over the end of $W$. Of course, the number
of ends must be zero modulo $2$.

The advantage of reducing to the case $b_1=1+b^+$ is that then, generically, 
all non-flat twisted reducibles in the moduli spaces are isolated.
Working directly with the original manifold $X$ would require dealing with
positive-dimensional families of twisted reducibles.
This technically more difficult situation has been studied by
Teleman \cite{Teleman1}. However,
it is not clear to this author whether one can expect stronger results
with such a direct approach.


After this paper was submitted the preprint \cite{nakamura1} appeared, which
addresses similar issues for closed $4$--manifolds, using
Seiberg--Witten theory.

{\bf Acknowledgement:} The author is grateful to the anonymous referee
for the careful reading of the manuscript and many suggestions.
He would also like to thank 
Ian Hambleton and Bj\o rn Jahren for helpful correspondence.

\section{Homology and cohomology with local coefficients}
\label{sec:cohomology-local}

This section contains mostly background material.

{\bf(I)} This part is concerned with singular (co)homology with local
coefficients. Let $X$ be any space.
For any bundle $E\to X$ of discrete Abelian groups we
denote by $C_*(X;E)$ the singular chain complex of $X$ with values in
$E$, as defined in \cite{Hatcher1}. A short exact sequence
\[0\to E'\to E\to E''\to0\]
of morphisms of such bundles induces a short exact sequence of chain
complexes
\[0\to C_*(X;E')\to C_*(X;E)\to C_*(X;E'')\to0\]
which in turn yields a long exact sequence relating the corresponding
homology groups $H_*(X;\cdot)$. Similar statements hold for the
singular cochain complexes and cohomology groups $H^*(X;\cdot)$.

Now let $p:\ti X\to X$
be any double covering
and $\ell\to X$ the associated  bundle
of infinite cyclic groups. Consider the $\z^2$--bundle
\[E:=\ti X\uset{\z/2}\times\z^2\]
over $X$, where $1\in\z/2$ acts on $\ti X$ by
flipping the sheets of the covering and on $\z^2$ by permuting the factors. Then
\[H_*(X;E\otimes G)=H_*(\ti X;G)\]
for any Abelian group $G$, and similarly for cohomology. There is a canonical
short exact sequence of bundles
\[0\to\ell\to E\to\z\to0\]
which induces a long exact sequence
\be{equation}\label{eqn:les-double-cover}
\cdots\to H_k(X;\ell)\to H_k(\ti X;\z)\oset{p_*}\to H_k(X;\z)\to
H_{k-1}(X;\ell)\to\cdots.
\end{equation}
We will use the notation $\lla$ (resp.\ $\dla$) for $\ell\otimes\R$
thought of as a real line bundle (resp.\ a bundle with discrete fibres)
over $X$. By the universal coefficients theorem
(see \cite[p.\,283]{Spanier}) one has
\[H^*(X;\ell)\otimes\R=H^*(X;\dla)\]
in each degree in which $H_*(X;\ell)$ is finitely
generated. There is a canonical isomorphism of bundles
$\R\oplus\dla\oset\approx\to E\otimes\R$, which induces an isomorphism
\[H^*(\ti X;\R)=H^*(X;\R)\oplus H^*(X;\dla).\]
The two summands correspond to the
$\pm1$ eigenspaces of the endomorphism of $H^*(\ti X;\R)$
induced by the involution of $\ti X$ (i.e.\ the action of $1\in\z/2$).
If $X$ is a smooth manifold then $H^*(X;\dla)$ can 
be computed as the de Rham cohomology associated to the flat bundle $\lla$
(see \cite{Bredon}). When working with de Rham cohomology it is 
natural to write $b_j(X;\lla)$ instead of $b_j(X;\ell)$, and similarly
for $b^+$.

There is also a relationship with mod~$2$ (co)homology, for arbitrary $X$:
The short exact sequence
\[0\to\ell\oset{\cdot2}\to\ell\to\z/2\to0\]
of bundles gives rise to a long exact sequence
\be{equation}\label{eqn:exact-cdot2}
\cdots\to H_q(X;\ell)\oset{\cdot2}\to H_q(X;\ell)\to H_q(X;\z/2)\to
H_{q-1}(X;\ell)\to\cdots
\end{equation}
as well as a similar sequence for cohomology. Furthermore,
because $\ell^*=\ell$ the universal coefficient theorem yields a split
short exact sequence
\be{equation}\label{eqn:univ-coeff}
0\to\text{Ext}(H_{q-1}(X;\ell),\z)\to H^q(X;\ell)
\to\text{Hom}(H_q(X;\ell),\z)\to0.
\end{equation}


\be{prop}\label{prop:order-two}
Let $X$ be any compact manifold (with or without boundary)
such that $H_1(X;\z/2)=\z/2$. Let $\ell\to X$ be any non-trivial bundle 
of infinite cyclic groups. Then $H_1(X;\ell)$ is a finite group of odd order,
hence by \Ref{eqn:univ-coeff} the group $H^2(X;\ell)$ contains no $2$--torsion.
\end{prop}

\proof Because $X$ is a manifold and $\ell$ is non-trivial, $H_0(X;\ell)=\z/2$.
Thus \Ref{eqn:exact-cdot2} yields an exact sequence
\[H_1(X;\ell)\oset{\cdot2}\to H_1(X;\ell)\to0.\]
Since $X$ is a {\em compact} manifold, $H_*(X;\ell)$ is finitely generated,
hence $H_1(X;\ell)$ must be a finite group on which multiplication by $2$
is an isomorphism.\square \\

We will now state a version of Poincar\'e duality for local
coefficients. Let $X$ be a closed topological $n$--manifold and $\cO_X\to X$
the orientation bundle, whose fibre over $x\in X$ is
\[\cO_x=H_n(X,X\setmin\{x\};\z).\]
Let $[X]\in H_n(X;\cO_X)$ be the fundamental class,
which is the unique class whose image in
\[H_n(X,X\setmin\{x\};\cO_X)=H_n(X,X\setmin\{x\};\cO_x)=\cO_x\otimes\cO_x=\z\]
is $1$ for every $x\in X$. Let $R$ be a commutative ring with identity.

\be{prop}\label{prop:poinc-duality}
For any closed topological
$n$--manifold $X$ and any bundle $E\to X$ of $R$--modules,
cap product with $[X]$ defines an isomorphism
\[H^p(X;E)\oset\approx\to H_{n-p}(X;E\otimes\cO_X)\]
for every $p$.
\end{prop}

\proof The proof in \cite{Hatcher1} for $R$--oriented $X$ and $E=R$
carries over with virtually no changes.\square

Other duality theorems for (co)homology with local coefficients can be
found in \cite{Spanier2}.

Now suppose $X$ is a closed oriented topological $n$--manifold and $\lla\to
X$ a bundle of infinite cyclic groups. Then it follows from
Proposition~\ref{prop:poinc-duality} and
the universal coefficient theorem \Ref{eqn:univ-coeff} (recalling that
$H_*(X;\ell)$ is finitely generated) that the 
intersection form $Q_{X,\ell}$ is unimodular. The same holds if $X$
has an integral homology sphere as boundary, as one can see by applying the
previous result to the double of $X$ and noting that the intersection
form of the double is the orthogonal sum of the intersection
forms of the two pieces.



{\bf(II)} In this part we use \v Cech cohomology.
Recall that for any paracompact space $X$ the first Chern class induces
an isomorphism between the group 
of isomorphism classes of complex line bundles over $X$ and
the cohomology group $H^2(X;\z)$. We will now give a
similar interpretation of $H^2(X;\ell)$. Let $\ti X,\lla$ be as in (I) and set
\[K:=\ti X\times_\zmt\co=\R\oplus\lla,\]
where $1\in\zmt$ acts on $\ti X$ by flipping the sheets and on $\co$ by
complex conjugation. Here $\co$ has the Euclidean topology, so that $K$ is a
real vector bundle over $X$.
Since conjugation is a field automorphism, $K$ is a 
bundle of fields isomorphic to $\co$. Let $K^*\subset K$ be the subspace
of non-zero vectors thought of as a bundle of multiplicative groups,
and let $\ck$ and $\ck^*$ denote the sheaves of
continuous sections of $K$ and $K^*$, resp.
By a {\em $K$--line bundle} we mean a bundle $L\to X$
such that each fibre $L_x$ is a $1$--dimensional vector space over $K_x$, 
and such that these data satisfy the usual axiom of local triviality.
A {\em local trivialization} of $L$ over an open subset $U\subset X$
is an isomorphism $L|_U\oset\approx\to K|_U$ of $K|_U$--modules.
An atlas of such local trivializations gives rise to a \v Cech cocycle
with values in $\ck^*$. Standard arguments show
that $L$ is classified up to isomorphism by the
corresponding cohomology class $\tc(L)\in H^1(X;\ck^*)$.
If $X$ is paracompact then the short exact sequence of sheaves
\be{equation}\label{eqn:exp-seq}
0\to\ell\to \ck\oset\exp\to \ck^*\to1
\end{equation}
yields an isomorphism $H^1(X;\ck^*)\oset\approx\to H^2(X;\ell)$, and we obtain:

\be{prop}\label{prop:twisted-c1}
For any paracompact space $X$
the characteristic class $\tc$ induces an isomorphism
between the group of isomorphism classes of $K$--line bundles
and the cohomology group $H^2(X;\ell)$.\square
\end{prop}

Note that $\La^2L=\lla$, so for the first Stiefel--Whitney
class one has
\[w_1(L)=w_1(\lla).\]
Furthermore, $\tc(L)$ maps to $w_2(L)$ under the homomorphism
$H^2(X;\ell)\to H^2(X;\z/2)$.

By a {\em Hermitian} $K$--line bundle we mean a $K$--line bundle
equipped with a Euclidean metric such that multiplication
with any unit vector in
$K_x$ is an orthogonal transformation of $L_x$, for any $x\in X$.

We now turn to the smooth category.
The proof of the following proposition is similar to that of
Proposition~\ref{prop:twisted-c1}.

\be{prop}\label{prop:smooth-twisted-c1}
For any smooth manifold $X$
the characteristic class $\tc$ induces an isomorphism
between the group of isomorphism classes of smooth Hermitian $K$--line bundles
and the cohomology group $H^2(X;\ell)$.\square
\end{prop}


Let $L\to X$ be a smooth Hermitian $K$--line bundle.
If $A$ is any (orthogonal) connection in $L$ then
its curvature $F_A$ is a $2$--form on $X$ with values in the bundle
$\so L$ of skew-symmetric endomorphisms of $L$. Under the
isomorphism $\lla\oset\approx\to\so L$ (defined by multiplication with elements
from $\lla$)
the closed form $F_A\in\Om^2(X;\lla)$ represents the image of
$-2\pi\tc(L)$ in $H^2(X;\dla)$.
(One can deduce this last statement
from the known case when $\ell$ is trivial by pulling
$A$ back to $\ti X$ and noting that $H^2(X;\dla)\to H^2(\ti X;\R)$ is
injective.)

\section{$\SO3$--connections and holonomy}
\label{sec:conn-hol}


Let $Z$ be a connected smooth $n$--manifold, possibly with boundary,
and $E\to Z$ an oriented, Euclidean
$3$--plane bundle. Fix $p>n$ and let $A$ be an (orthogonal) $L^p_{1,\text{loc}}$
connection in $E$. Let $\Ga_A$ denote the group of $L^p_{2,\text{loc}}$
automorphisms of $E$ which preserve $A$. Just as for smooth connections,
$\Ga_A$ is isomorphic to the centralizer of the holonomy group
$\holza\subset\aut(E_z)\approx\SO3$ at any point $z\in Z$.
Recall that any positive-dimensional proper closed subgroup
of $\SO3$ is conjugate to either $\U1$
or $\O2$, and these subgroups have centralizer $\U1$ and $\z/2$, resp.
We will call the connection $A$
\be{itemize}
\item {\em irreducible} if $\Ga_A=\{1\}$, otherwise {\em reducible},
\item {\em Abelian} if $\Ga_A\approx\U1$,
\item {\em twisted reducible} if $\Ga_A\approx\z/2$.

\end{itemize}

Now suppose $A$ is smooth. Then $A$ is reducible if and only if it
preserves a rank~$1$ subbundle $\lla\subset E$. If in addition
$A$ is not flat then $\lla$ is unique (because a non-flat connection $A$
has holonomy 
close to but different from $1$ around suitable small loops in $Z$).
In that case $A$ is 
Abelian if $\lla$ is trivial and twisted reducible otherwise.

Now suppose $Z$ is compact. Let $\ca$ denote the affine Banach space
consisting of all $L^p_1$ connections in $E$ and let $\cg$ be the 
Banach Lie group of all $L^p_2$ automorphisms of $E$.
Then $\cg$ acts smoothly on $\ca$ and we denote the quotient space by
$\cb=\cb(E)$. It follows easily from the local slice theorem (see
\cite[p.\ 132 and p.\,192]{DK} and \cite[Section~2.5]{Fr13}) that
$\cb$ is a regular topological space. Since $\cb$ is also second
countable, it is metrizable by the Urysohn metrization theorem \cite{Kelley}.
Hence $\cb$ is paracompact, and the same holds for any subspace of $\cb$.

Let $\ca^*\subset\ca$ be the subset of irreducible connections.
Then $\cb^*:=\ca^*/\cg$ is a Banach manifold.
In the proof of the theorem we will take $p$ to be an even integer, to make
sure that $\cb^*$ possesses smooth partitions of unity. (In \cite{Lang1}
the existence of smooth partitions of unity is established for
paracompact Hilbert manifolds. The proof carries over to paracompact
Banach manifolds $B$ modelled on a Banach space $(E,\|\cdot\|)$ such that
$\|\cdot\|^t$ is a smooth function on $E$ for some $t>0$. This
includes $B=\cb^*$ when $p$ is an even integer, with $t=p$.)

Recall that the Lie group $\aut(E_z)\approx\SO3$ has
a non-trivial double covering
\be{equation}\label{eqn:taut}
\waut\to\aut(E_z),
\end{equation}
where $\waut$ is isomorphic to the group $\Sp1$ of unit quaternions.
Let $\cg$ act on $\aut(E_z)$ by conjugation with $u(z)$ and on $\waut$ by
conjugation with any lift of $u(z)$ to $\waut$. Then the covering map
\Ref{eqn:taut} is $\cg$--equivariant.
It follows from the local slice theorem that $\ca^*\to\cb^*$ is a 
principal $\cg$--bundle, hence
\be{equation}\label{eqn:double-cover}
\ca^*\times_\cg\waut\to\ca^*\times_\cg\aut(E_z)
\end{equation}
is a double covering. Now let $\ga:S^1\to Z$ be a loop based at $z$.
Pulling back
\Ref{eqn:double-cover} by the smooth map
\[\cb^*\to\ca^*\times_\cg\aut(E_z),\quad[A]\mapsto[A,\hola]\]
yields a double covering of $\cb^*$.
We will now show that this
extends to a double covering $\xg\to U_\ga$, where $U_\ga\subset\cb$
contains $\cb^*$ as well as certain reducibles.

\be{defn}
\be{description}
\item[(i)]Let $U_\ga\subset\cb$ be the subspace
consisting of those $[A]$ such that there are two points in 
$\ca\times_\cg\waut$ lying above $[A,\hola]\in\ca\times_\cg\aut(E_z)$.
\item[(ii)]Let $\xg\subset\ca\times_{\cg}\waut$
be the subspace consisting of those $[A,g]$ such that $[A]\in U_\ga$ 
and $g\in\waut$ is a lift of $\hola$.
\end{description}
\end{defn}

{\em Remark:} Note that $[A]\in\cb$ lies in the complement of $U_\ga$ if and only if
there exists a $u\in\Ga_A$ such that $u$ interchanges the two points
in $\waut$ lying above $\hola$, or equivalently, such that $u(z)$ and $\hola$
are both reflections and have perpendicular axes of rotation.

\be{prop}\label{prop:red-in-Uga}
Let $[A]\in\cb$.
\be{description}
\item[(i)]If $A$ is Abelian, then $[A]\in U_\ga$.
\item[(ii)]Let $A$ be twisted reducible and let $\lla\subset E$ be
the $1$--eigenspaces
of the non-trivial element of $\Ga_A$. Then $[A]\in U_\ga$ if and only if
$\ga^*\lla$ is trivial.
\end{description}
\end{prop}

Note that elements of $\cg$ are of class $C^1$ by the Sobolev embedding
theorem, hence the subbundle $\lla\subset E$ in (ii) is of class $C^1$.

\proof (i) If $\Ga_A\approx\U1$ then $\Ga_A$ is the centralizer of any 
non-trivial element $x\in\Ga_A$ with $x^2\neq1$. Hence $\hol_z(A)\subset\Ga_A$, so $[A]\in U_\ga$ by the above remark.

(ii) Since $A$ preserves the subbundle $\lla$, the holonomy $\hol_\ga(A)$
acts as $\eps=\pm1$ on the fibre $\lla_z$. Therefore, $[A]\in U_\ga$ if and
only if $\eps=1$, or equivalently, if $\ga^*\lla$ is trivial.\square

\be{prop}\label{prop:double-covering}
$U_\ga$ is an open subset of $\cb$, and
the canonical projection $\xg\to U_\ga$ is a double covering.
\end{prop}

\proof We give a proof which does not require the local
slice theorem.
After choosing a framing of $E_z$ we can identify the
covering \Ref{eqn:taut}
with the adjoint representation $\Sp1\to\SO3$.
Fix $A\in\ca$ with $[A]\in U_\ga$ and a lift $q\in\Sp1$ of $\hola$.
For $\eps>0$ set
\[P_\eps:=A+\oset\circ D_\eps,\]
where $\oset\circ D_\eps\subset L^p_1(Z;\so E)$ is the open $\eps$--ball
about the origin. Let $\pi:\ca\to\cb$ be the projection. This is an
open map, since $\cb$ is the quotient of $\ca$ with respect to a group action.
Hence $\pi(P_\eps)\subset\cb$ is an open neighbourhood of $[A]$.
If $B\in P_\eps$ with $\eps$ sufficiently small then
$\hol_\ga(B)\cdot\hola\inv$ will not be a reflection and so has a unique
lift $g(B)\in\Sp1$ with positive real part. Then
\[f(B):=g(B)q\]
is a lift of $\hol_\ga(B)$. A simple convergence argument shows that
if $\eps$ is sufficiently small and $B\in P_\eps$, $u\in\cg$ are such that
$u(B)\in P_\eps$ then
\[f(u(B))=u\cdot f(B).\]
For such $\eps$ we have $\pi(P_\eps)\subset U_\ga$, and 
the map $[B]\mapsto[B,f(B)]$ is a continuous section
of $\Xi_\ga$ over $\pi(P_\eps)$.
Changing the sign of $f$ yields a different
section and altogether a trivialization
of $\Xi_\ga$ over $\pi(P_\eps)$.\square

\section{Moduli spaces and twisted reducibles}
\label{sec:tw-red}

Let $W$ be any oriented, connected, Riemannian $4$--manifold
with one cylindrical
end $\rpy$, where $Y$ is an integral homology sphere. (Thus, the complement
of $\rpy$ is compact). Let $E\to W$
be an oriented Euclidean $3$--plane bundle. Choose a trivialization
of $E|_{\rpy}$. For any non-degenerate
flat connection $\rho$ in the product $\SO3$--bundle $E_0\to Y$
let $M(E,\rho)$ denote the moduli space
of instantons in $E$ that are asymptotic to $\rho$ over the end. We briefly
recall the construction of this moduli space, following \cite{D5,Fr13}.
Choose a smooth reference connection $\aref$ in $E$ whose restriction to
the $\rpy$ is the pull-back of $\rho$. Introduce the space
\[\ca=\ca(E,\rho):=\aref+L^{p,w}_1(W;\so E)\]
of Sobolev connections, where $w$ is a small exponential weight
as in \cite[Subsection~2.1]{Fr13} (which is
actually only needed when $\rho$ is reducible). There is a Banach Lie group
$\cg$ (consisting of certain $\lloc p2$ gauge transformations)
acting on $\ca$, and $M(E,\rho)$ is the subspace of the quotient space
$\cb:=\ca/\cg$ consisting of all $[A]$ satisfying $F^+_A=0$.

If $u:Y\to\SO3$ then the moduli spaces with limits
$\rho$ and $u(\rho)$, resp.,
can be identified if $u$ is null-homotopic; otherwise
the expected dimensions of these moduli spaces differ by $4\deg(u)$.
Let $\fl_Y$ denote the space of gauge equivalence classes of flat
connections in $E_0$, and let $\fl^*_Y$ be the irreducible part of $\fl_Y$.
It will be convenient to denote a moduli space $M(E,\rho)$ of expected dimension
$d$ by $M_{\al,d}$, where $\al=[\rho]\in\fl_Y$.
In the particular case when $\rho$ is trivial, however, we will usually
label the
moduli space by $k=-p_1(E,\rho)\in H^4_c(W;\z)=\z$, where $p_1(E,\rho)$ is
the relative Pontryagin class. Note that as $\rho$ varies, $k$
runs through a set of the form $k_0+4\z$, $k_0\in\z$.
Thus, $M_k$ will denote the moduli space
with trivial limit and expected dimension
\[\dim M_k=2k-3\del(W),\]
where
\be{equation}\label{eqn:deldef}
\del(W):=1-b_1(W)+b^+(W).
\end{equation}
If $M_k$ is non-empty then for every $[A]\in M_k$ one has
\be{equation}\label{eqn:fafa}
8\pi^2k=\int_W\tr(F_A\wedge F_A)=\int_W|F^-_A|^2\ge0.
\end{equation}

After perturbing the Riemannian metric on $W$ in a small ball we may assume that
there is no 
$[A]\in M_k$, for any $k>0$, such that $A$ preserves
a real line bundle $\lla\subset E$ with $b^+(W;\lla)>0$.
(This can be proved along the
same lines as the untwisted case \cite[Corollary~4.3.15]{DK},
cf.\ \cite[Lemma~2.4]{KM3}.)

For the remainder of this section assume
\be{equation}\label{eqn:k0delW0}
k>0,\qquad \del(W)=0.
\end{equation}
Then $M_k=M_{\theta,2k}$, where $\theta\in\fl_Y$ is the class of
trivial connections. Let $M_k^*,M_k\red,M_k\tred$ be the subsets of $M_k$
consisting of the irreducible, reducible, and twisted reducible points,
resp.

\be{prop}\label{prop:twisted-red-class}
There is a canonical bijection between $M\tred_k$ and the set $P$ of equivalence
classes of pairs $(\ell,c)$, where $\ell\to W$ is a non-trivial
bundle of infinite cyclic
groups, $c\in H^2(W;\ell)$, and such that for $\lla:=\ell\otimes\R$ one has
\[b^+(W;\ell)=0,\quad w_1(\lla)^2+[c]_2=w_2(E),\quad c^2=-k,\]
where $[c]_2$ denotes the image of $c$ in $H^2(W;\z/2)$.
\end{prop}

Here two such pairs $(\ell,c)$, $(\ell',c')$ are deemed equivalent if
there is an isomorphism 
$\ell\oset\approx\to\ell'$
such that $c\mapsto c'$ under the 
induced isomorphism $H^2(W;\ell)\oset\approx\to H^2(W;\ell')$.

\proof (i) To define this bijection, let $[A]\in M\tred_k$. 
We may assume $A$ is smooth. Since $A$ is not flat, it preserves
a unique non-trival rank~$1$ subbundle $\lla\subset E$. Let
$K=\R\oplus\lla$ be the corresponding bundle of fields as in
Section~\ref{sec:cohomology-local}. The orthogonal 
complement $L\subset E$ of $\lla$ is in a canonical way a $K$--line bundle. The
module structure is given as follows: For $x\in W$, $(a,b)\in\R\oplus\lla_x$,
$v\in L_x$ set
\be{equation}\label{eqn:K-module}
(a,b)\cdot v:=av+b\times v,
\end{equation}
where $b\times v$ is the cross product in the $3$--dimensional, oriented,
Euclidean vector space $E_x$.
Let $\ell\subset\lla$ denote the lattice of vectors of integer length
and set $c:=\tc(L)\in H^2(W;\ell)$. It is clear that different
representatives $A$ of the same point in $M\tred_k$ are mapped to equivalent
pairs $(\ell,c)$.

We now verify that $(\ell,c)$ has the required properties. 
By choice of metric on $W$ we must have $b^+(W;\ell)=0$. Furthermore,
\[w_2(E)=w_2(\lla\oplus L)=w_1(\lla)\cup w_1(L)+w_2(L)=w_1(\lla)^2+[c]_2.\]
Secondly, let $B$ denote the connection
in $L$ induced by $A$. Then $F_B$ takes values
in $\lla$, and one has
\[\tr(F_A\wedge F_A)=-2F_B\wedge F_B\in\Om^4(W).\]
Since $F_B$
decays exponentially, we obtain
\[\int_W\tr(F_A\wedge F_A)=-2\int_WF_B\wedge F_B=-8\pi^2c^2,\]
hence $c^2=-k$.

(ii) Now suppose $A,A'\in\ca$ are smooth connections
representing points in $M\tred_k$,
and that the corresponding
pairs $(\ell,c)$, $(\ell',c')$ are equivalent through an isomorphism
$f:\ell\oset\approx\to\ell'$. Let $E=\lla\oplus L$ and
$E=\lla'\oplus L'$ be the splittings preserved by $A$ and $A'$, resp.,
and let $K,K'$ be the bundles of fields corresponding to $\lla,\lla'$, resp.
Let $\phi:K\to K'$ be the isomorphism induced by $f$.
By means of $\phi$, we turn $L'$ into an Hermitian $K$--line bundle
which we denote by $L'_\phi$. It is easy to check that 
$f_*(\tc(L'_\phi))=\tc(L')$, so by
Proposition~\ref{prop:smooth-twisted-c1} there is an isomorphism
$\psi:L\to L'_\phi$ of Hermitian $K$--line bundles.
Combining $\phi|_\lla$ and $\psi$ we obtain an isomorphism of Euclidean
vector bundles
\[u:E=\lla\oplus L\to\lla'\oplus L'=E.\]
To see that $u$ preserves orientations, let $a\in\lla_x$ and $b\in L_x$ be
of unit length. Then $(a,b,a\times b)$ is a positive orthonormal basis for
$E_x$. Under $u$ this is mapped to $(\phi(a),\psi(b),\phi(a)\times\psi(b))$,
which is also a positive orthonormal basis.

We may assume $A$ and $A'$ are in temporal gauge. Then $L$ and $L'$
will be translationary invariant over the end $W^+:=\R_+\times Y$
with respect to the chosen trivialization of $E|_{W^+}$. Let $v$ be
the non-trivial element of $\Ga_A$ and let $A|_{W^+}=d+a$, where $d$
denotes the product connection. Then over the end one has
\[0=d_Av=dv+av-va.\]
Since $v$ is translationary invariant over the end and
\[\int_{[t,t+1]\times Y}|a|^p\to0\quad\text{as $t\to\infty$,}\]
we conclude that $dv=0$ on $W^+$. The same holds for the non-trivial
element of $\Ga_{A'}$.
Hence
\[L|_{W^+}=W^+\times C,\qquad L'|_{W^+}=W^+\times C'\]
for some $2$--dimensional subspaces $C,C'\subset\R^3$. Now
$\psi|_{W^+}$ is given by a smooth map
\[\ti\psi:W^+\to\text{SO}(C,C'),\]
where $\text{SO}(C,C')\approx S^1$ is one specific component of the space of
linear isometries $C\to C'$, the component being determined by the
isomorphism $f:\ell\to\ell'$. But every map $C\to S^1$ is
null-homotopic, since $H^1(W^+;\z)=0$. We may therefore choose the
isomorphism $\psi$ such that $\ti\psi$ is constant on
$[1,\infty)\times Y$, say. Then $du=0$ on $[1,\infty)\times Y$, so
$u\in\cg$. Set $A'':=u\inv(A')\in\ca$.

Recall that the cross product on $E$ defines a canonical isomorphism
$E\oset\approx\to\so E$. Under this isomorphism, the difference
$b:=A''-A$ is a $1$--form with values in $\lla$. More precisely,
$b\in\llw pw1(W;\La^1\otimes\lla)$. Moreover,
\[F_{A''}=F_A+db,\]
so $d^+b=0$. Since $b_1(W;\ell)=0$ by \Ref{eqn:del-indep}
there is a section $\xi\in\llw pw2(W;\lla)$
such that $d\xi=b$. Set $v=\exp(\xi)$. Then $v(A'')=A$, so $A$ and $A'$
represent the same point in $M\tred_k$.

(iii) We will now show that every class $[\ell,c]\in P$ is the image of some
point $[A]\in M\tred_k$. Define $\lla,K$ in terms of $\ell$ as 
in Section~\ref{sec:cohomology-local}.
Choose a $K$--line bundle $L\to W$ with $\tc(L)=c$. The hypotheses on $\ell,c$
imply that $\lla\oplus L$ and $E$ have the same second Stiefel--Whitney class,
hence these bundles are isomorphic (see \cite[p.\,674]{Dold-Whitney} and
\cite[Theorem E.8]{FU}); we will identify them. 
Since $L$ is trivial
over the end of $W$, there is an orthogonal connection $A'$ in $E$
which respects the given splitting and is flat over the end of $W$.
Since $d^+:\Om^1(W;\lla)\to\Om^+(W;\lla)$ induces a surjective map
$L^{p,w}_1\to L^{p,w}$ between Sobolev spaces with a small positive weight
(cf.\ the proof of \cite[Prop.\ 5.1.2]{Fr13}),
there is an $a\in L^{p,w}_1$ such that $A:=A'+a$ satisfies
$F_A^+=0$. Clearly, $[A]\in M\tred_k$ is mapped to $[\ell,c]$.\square
\\

Now fix $\ell\to W$ and let $P_\ell$ be the set of points in $P$ of the
form $[\ell,c]$. Suppose $P_\ell\neq\emptyset$ (which implies
$b^+(W;\ell)=0$) and choose a $c$ with $[\ell,c]\in P_\ell$.
Let $\ctl$ be the torsion subgroup of $H^2(W;\ell)$ and 
for any $v\in H^2(W;\ell)$ let $\bar v$ denote the image of $v$ in 
$H^2(W;\ell)/\ctl$. Set
\[P_c:=\{\{r,s\}\subset H^2(W;\ell)/\ctl\st r\cdot s=0;\;r+s=\bar c\},\]
where $\{r,s\}$ means the {\em unordered} set.

\be{prop}\label{prop:cardinality}
$|P_\ell|=|2\ctl|\cdot|P_c|$.
\end{prop}

Here $|\cdot|$ denotes the cardinality of the given set. Note that $2\ctl$
has even order if and only if $H^2(W;\ell)$ contains an element of order $4$.

\proof Let $\tpl$ be the set of all $v\in H^2(W;\ell)$ such that
$[\ell,v]\in P$. Set
\[\al:\tpl\to P_\ell,\quad v\mapsto[\ell,v].\]
Since the only non-trivial automorphism of $\ell$ is given by multiplication
by $-1$, we have
\[\al(v)=\al(v')\iff v=\pm v'.\]
Because $k\neq0$ it follows that $\al$ is two-to-one, hence
\[|\tpl|=2|P_\ell|.\]
Now let $\tpc$ be the set of all ordered pairs $(r,s)$ such that
$\{r,s\}\in P_c$. Because $k\neq0$ one has $r\neq s$ for all such $r,s$,
hence
\[|\tpc|=2|P_c|.\]
It follows from the long exact sequence
\be{equation}\label{eqn:coh-bockstein}
\cdots\to H^2(W;\ell)\oset{\cdot2}\to H^2(W;\ell)\to H^2(W;\z/2)\to\cdots
\end{equation}
(see Section~\ref{sec:cohomology-local}) that the map
\[\tpl\to\tpc,\quad v\mapsto
\left(\frac{\bar c+\bar v}2,\frac{\bar c-\bar v}2\right)\]
induces a bijection $\tpl/2\ctl\to\tpc$, where $2\ct_\ell$ acts on $\tpl$
by translation, hence 
\[|\tpl|=|2\ctl|\cdot|\tpc|\]
and the proposition is proved.\square \\

\section{Local structure around twisted reducibles}
\label{sec:local-str-twr}

We continue the discussion of the previous section, under the assumptions
\Ref{eqn:k0delW0}.

We do not know if the twisted reducibles in $M_k$ are regular points of $M_k$
for a generic tubular end metric on $W$ (although there is a generic
metric theorem of this kind for closed $4$--manifolds,
see \cite[Lemma~2.4]{KM3}). However,
regularity of these reducibles can be achieved by a simple local perturbation
of the instanton equation which is similar in spirit to 
that used in \cite[p.\,292]{D1}. To describe this perturbation, 
let $M_k\subset\cb=\ca/\cg$ as in Section~\ref{sec:tw-red}, and suppose
${B}\in\ca$ satisfies $F^+_{B}=0$ and preserves a splitting $E=\lla\oplus L$,
where $\lla$ is a non-trivial real line bundle. 
Then the non-trivial element of $\Ga_B$ acts on any fibre of
$\lla\oplus L$ by $(a,b)\mapsto(a,-b)$.
For any $\eps>0$ set
\[S_{0,\eps}=\{a\in\llw pw1(W;\so E)\st d^*_{B}a=0,\;\|a\|_{\llw pw1}<\eps\},\]
where the Sobolev norm is defined in terms of ${B}$. This norm
is equivalent to the corresponding
norm defined by the reference connection $\aref$ because
of the Sobolev embedding $L^p_1\subset L^\infty$ in $\rf$. (Recall that we
are assuming $p>4$.) If $\eps$ is sufficiently small then
$S_\eps:={B}+S_{0,\eps}$
is a local slice to the action of $\cg$. This means, firstly, that
there is an open neighbourhood $U$ of $1\in\cg$ such that
\[U\times S_\eps\to\ca,\quad(u,A)\mapsto u(A)\]
is a diffeomorphism onto an open subset of $\ca$, and secondly, that the
projection $S_\eps/\Ga_{B}\to\cb$ is injective. Then $S_\eps/\Ga_{B}$ maps
homeomorphically onto an open neighbourhood of $[{B}]$ in $\cb$, and 
the irreducible part of $S_\eps/\Ga_{B}$ maps diffeomorphically onto an open
subset of $\cb^*$.
The operator
\[-d_{B}^*+d^+_{B}:\Om^1\to\Om^0\oplus\Om^+,\]
acting on forms on $W$ with values in $\so E\approx E$,
induces a Fredholm operator
$\cD:\llw pw1\to L^{p,w}$ whose index is the expected dimension of $M_k$,
i.e.\ $\ind(\cD)=2k>0$. Therefore, there is a compact operator $P$
such that $\cD+P$ is surjective. We will choose such a $P$ of a
particular kind. To describe this, first note that
\be{equation}\label{eqn:cDllaL}
\cD=\cD_\lla\oplus \cD_L,
\end{equation}
where $\cD_\lla$ and $\cD_L$ act on forms with values in $\lla$ and $L$,
resp. Now, $\cD_\lla$ is an isomorphism, because $\lla$ is
non-trivial and $b^+(W;\lla)=\del(W)=0$. Therefore, $\cD+P$ is surjective if
$P$ is given by
\[Pa=\sum_{j=1}^r\la a,\phi_j\ra_{L^2}\cdot\om_j,\]
where $r$ is the dimension of the cokernel of $\cD_L$, and
$\phi_j\in\Om^1_c(W;L)$, $\om_j\in\Om^+_c(W;L)$ are suitably chosen.
Choose a smooth function $\ka:[0,\infty)\to[0,\infty)$ such that $\ka(t)=1$
for $t\le\eps/3$ and $\ka(t)=0$ for $t\ge2\eps/3$. For any $a\in S_{0,\eps}$
set
\[\frakp({B}+a):=\ka(\|a\|_{\llw pw1})\cdot Pa.\]
Then $\frakp$ is a smooth $\Ga_{B}$--equivariant map
$S_\eps\to\Om^+_c(W;\so E)$.
Moreover, $\frakp$ extends uniquely to a smooth $\cg$--equivariant
map $\ca\to\lw pw(W;\La^+\otimes\so E)$ which vanishes outside $\cg S_\eps$.
This extension will also be denoted $\frakp$.

The perturbed instanton equation that we have in mind is then
\be{equation}\label{eqn:inst-pert}
F^+_{A}+\frakp(A)=0,
\end{equation}
for $A\in\ca$.
Clearly, the linearization of this equation at $B$ is surjective, since
it restricts to $d^+_B+P$ on $\ker d^*_B$.
Note that adding the perturbation $\frakp$
does not affect the compactness properties of the corresponding moduli
space. If we take $\eps>0$ sufficiently small, then the classification
of twisted reducibles in Proposition~\ref{prop:twisted-red-class}
is also not affected.

More generally, we may add one such local perturbation $\frakp$ for
each of a finite number of twisted reducibles in $\cb$.
Usually, the perturbations will be suppressed from notation.

Having resolved the regularity issue, we now describe the local
structure around a regular twisted reducible in $M_k$.

In the next lemma $Z$ will denote a compact, connected
codimension~$0$ submanifold
of $W$. Consider the double covering $\Xi_\ga\to U_\ga$ associated to the bundle
$E|_Z$ and a loop $\ga:S^1\to Z$ based at $z\in Z$.

\be{lemma}\label{lemma:xi-link}
Suppose $[B]$ is a regular point of $M_k$ such that $B$ preserves a non-trivial
real line bundle $\lla\subset E$. Then under the restriction
map $R:M_k^*\to U_\ga$, the pull-back of the double covering $\Xi_\ga\to U_\ga$
is trivial over the link of $[B]$ in $M_k$ if and only if
$\ga^*\lla$ is trivial.
\end{lemma}
The fact that $\Xi_\ga\to U_\ga$ is a double covering was
proved in Proposition~\ref{prop:double-covering}.
By the ``link'' we mean the boundary $\prtl N\approx\mathbb{RP}^{2k-1}$
of a compact neighbourhood $N$ of $[B]$ in $M_k$ to be constructed in 
the proof.

\proof If $\ga^*\lla$ is trivial then $[B]\in U_\ga$ by
Proposition~\ref{prop:red-in-Uga}, so there is a well-defined restriction
map
\[\bar R:M_k^*\cup\{[B]\}\to U_\ga.\]
Since $\Xi_\ga\to U_\ga$ is locally trivial, $\bar R^*\Xi_\ga$ is trivial
on a neighbourhood of $[B]$.

Now suppose $\ga^*\lla$ is non-trivial. 
We may assume $B\in\ca$ is smooth. Recall that the kernel $K$
of the operator \Ref{eqn:cDllaL} consists entirely of forms with values in
$L$. Therefore the non-trivial element of $\Ga_B$ acts as $-1$ on $K$. 

For a small $r>0$ let $D_r\subset K$ be the closed $r$--ball around the
origin with respect to some $\Ga_B$--invariant inner product on $K$. By the
local slice theorem there is a smooth $\Ga_B$--equivariant embedding
\[Q:=B+D_r\to\ca\]
whose composition with the projection $\ca\to\cb$
induces a homeomorphism of $Q/\Ga_B$ onto a compact neighbourhood
$N$ of $[B]$ in $M_k$.

Let $\ti Q\to Q$ be the pull-back of the double covering \Ref{eqn:taut} under
$\hol_\ga:Q\to\aut(E_z)$. Since $Q$ is contractible, $\ti Q\to Q$ is a trivial
double covering. There is now a commutative diagram
\be{equation*}
\be{array}{ccc}
\prtl\ti Q/\Ga_B & \to & R^*\Xi_\ga\\
\downarrow & & \downarrow\\
\prtl Q/\Ga_B & \to & M_k^*
\end{array}
\end{equation*}
where the horizontal maps are the embeddings induced by $Q\to\ca$.
The image of the
bottom map is $\prtl N$, so what we need to show is that the left-most
map is a non-trivial covering. Since $\ga^*\lla$ is non-trivial,
$h:=\hol_\ga(B)$ acts as $-1$ on $\lla_z$ and hence by a reflection on $L_z$.
In a suitable orthogonal basis for $E_z$
the two lifts of $h\in\SO3$ to $\Sp1$ are $\pm j$,
and $\si$ acts on $\Sp1$ by conjugation with $i$.
Since $iji\inv=-j$, we see that $\si$ interchanges the two points in 
$\ti Q$ lying above $B$. Thus we can identify $\ti Q\to Q$ with
\[D^{2k}\times\{\pm1\}\to D^{2k},\]
where $D^{2k}$ is the unit disk in $\R^{2k}$, and $\si$ acts
on $D^{2k}\times\{\pm1\}$ by $(x,t)\mapsto(-x,-t)$. Restricting to
$\prtl D^{2k}=S^{2k-1}$ and dividing out by $\Ga_B$ we obtain the usual covering
$S^{2k-1}\to\mathbb{RP}^{2k-1}$, which is non-trivial.\square

\section{Abelian instantons}
\label{sec:abelian-inst}

Let $M_k$ be as in Section~\ref{sec:tw-red}, assuming
\Ref{eqn:k0delW0}. In this section we use the unperturbed instanton equation.

We will need to deal with both reducibles in $M_k$ and reducibles
appearing in weak limits of sequences in $M_k$. Such reducibles are
contained in the set
\[M\red:=\coprod_{s\ge0}M\red_{k-4s}.\]
Let $M\ared$ and $M\tred$ be the subsets of $M\red$ consisting of
Abelian and twisted reducibles, resp.
If $b^+(W)>0$ then $\mared$ is empty, by choice of the metric on $W$.
If $b^+(W)=0$ (in which case $b_1(W)=1$)
then $\mared$ consists of a finite disjoint union of circles. Let $Z\subset W$
be any compact, connected, codimension~$0$ submanifold and $\ga$ a loop in $Z$.
The restriction map
\[R:\mared\to U_\ga\subset\cb(E|_Z)\]
maps each circle $S\subset\mared$ onto either a circle or a point, depending
on whether $H^1(W;\R)\to H^1(Z;\R)$ is non-zero or not. Moreover,
the double covering $\Xi_\ga\to U_\ga$
pulls back to a trivial covering of any $S$
if and only if the class in $H_1(W;\z)/\torsion$ represented by $\ga$ 
is divisible by $2$; however, we will make no use of this.
Note that different circles $S$ have disjoint images in $\cb(E|_Z)$, by the 
unique continuation property of self-dual closed $2$--forms (applied to the
curvature forms).

\section{Three lemmas on Banach manifolds}
\label{sec:banach}

The results of this section will be used in
Section~\ref{subsec:choosing-sec} below.

\be{lemma}\label{lemma:cont-smooth}
Let $B$ be a smooth (Hausdorff)
Banach manifold which admits smooth partitions of unity.
Let $L\to B$ be a smooth real $2$--plane
bundle. If $L$ admits a continuous non-vanishing section, then it 
admits a smooth non-vanishing section.
\end{lemma}

\proof Of course, this is well known if $B$ is finite-dimensional (and at
least in that case it holds for bundles of any finite rank).
For general $B$ one can use the following \v Cech cohomology argument:

Choose a smooth Euclidean metric on $L$. Set $\lla:=\La^2L$ and let
$K:=\R\oplus\lla$ be the associated bundle of fields as in
Section~\ref{sec:cohomology-local}.
Then $L$ has a canonical structure as a smooth $K$--line bundle as defined
in \Ref{eqn:K-module}. Clearly, $L$ is trivial as a smooth (resp.\ continuous)
$K$--line bundle (meaning $L\approx K$) if and only if $L$ admits a
non-vanishing smooth (resp.\ continuous) section.

Let $\ck^*$ and $\ck^*_\infty$ denote the sheaves of continuous and smooth
sections of $K^*$, resp. Then $L$ is determined up to smooth isomorphism
by its class in $H^1(B;\ck^*_\infty)$. But the inclusion $\ck^*_\infty\to\ck^*$
induces an isomorphism
\[H^1(B;\ck^*_\infty)\oset\approx\to H^1(B;\ck^*),\]
as is easily seen by considering the morphism between the exponential
short exact sequences for $\ck^*_\infty$ and $\ck^*$ (see \Ref{eqn:exp-seq}).
\square \\

In the following two lemmas, $B$ will be a metric space and $R$ a
compact subspace. The open subspace
$B^*:=B\setminus R$ of $B$ will have the structure of a smooth Banach manifold 
admitting smooth partitions of unity.

\be{lemma}\label{lemma:non-van1}
Suppose $R$ is a finite set. Let $\Theta\to B$ be a Euclidean (real)
line bundle. Let
$\Theta|_{B^*}$ have the obvious smooth structure. Then there exists a
smooth section of $\Theta|_{B^*}$ which is nowhere zero in
$B^*\cap V$ for some neighbourhood $V$ of $R$ in $B$.
\end{lemma}

\proof Choose an open neighbourhood $N$ of $R$ and a section $\si$ of
$\Theta|_N$ which has (pointwise) unit length. Then $\si$ is smooth in
$N\cap B^*$.
Since $B$ is a normal space there are disjoint open
neighbourhoods $V,V'$ of $R$ and $B\setminus N$, resp. Using a
smooth partition of unity of $B^*$ subordinate to the open cover
$\{V',B^*\cap N\}$ one can construct a smooth section $s$ of $\Theta$
over $B^*$ which agrees with $\si$ on $B^*\cap V$.
In particular, $s$ is nowhere zero in $B^*\cap V$.\square

\be{lemma}\label{lemma:non-van2}
Suppose $R$ is the disjoint union of three sets,
\[R=R_0\sqcup R_1\sqcup R_2,\]
where $R_1$ and $R_2$ are finite sets and $R_0$ is a finite disjoint union of
subspaces each of which is homeomorphic to a circle. For $i=1,2$ let
$\Theta_i\to B^*\cup R_0\cup R_i$ be a Euclidean line bundle.
Let $\hth$ be the direct sum
of the restrictions of $\Theta_1$ and $\Theta_2$ to $B^*\cup R_0$. Then there
exists a smooth section of $\hth$ over $B^*$
which is nowhere zero in $B^*\cap V$
for some neighbourhood $V$ of $R$ in $B$.
\end{lemma}

\proof Choose pairwise disjoint open sets $H_0,H_1,H_2$ in $B$ such
that $R_i\subset H_i$ for $i=0,1,2$. It is easy to see that
$\hth|_{R_0}$ admits a non-vanishing section $\si$. Since $R_0$ is
compact we can cover $R_0$ by finitely may open sets $U_j$ in $H_0$ such
that $\hth|_{U_j}$ is trivial for each $j$. By the Tietze extension
theorem there is a section $\si_j$ of $\hth|_{U_j}$ which agrees with
$\si$ on $R_0\cap U_j$. Patching
together the sections $\si_j$ by means of a partition of unity yields
a section $\ti\si$ of $\Theta$ over
$U:=\cup_jU_j$ such that $\ti\si=\si$ on $R_0$. Then the locus $N_0$
where $\ti\si\neq0$ is an open neighbourhood of $R_0$ in $B$.
By Lemma~\ref{lemma:cont-smooth} there exists a non-vanishing smooth
section $s_0$ of $\hth$ over $B^*\cap N_0$.

For $i=1,2$ choose a unit length section $\tau_i$ of $\Theta_i$ over some
open neighbourhood $N_i$ of $R_i$ in $H_i$. Combining $\tau_i$ with the
zero-section of $\Theta_{3-i}$ yields a smooth non-vanishing section
$s_i$ of $\hth$ over $B^*\cap N_i$.

Set $N:=N_0\cup N_1\cup N_2$, which is an open neighbourhood of $R$ in
$B$. By means of a smooth partion of unity as in the proof of
Lemma~\ref{lemma:non-van1} we can then construct a smooth section $s$
of $\hth|_{B^*}$ which agrees with $s_i$ in $B^*\cap V\cap N_j$ for some
neighbourhood $V$ of $R$ in $N$. In particular, $s$ is nowhere zero in
$B^*\cap V$.\square

\section{Proof of theorem}
\label{sec:proof-thm}

Assuming the hypotheses of the
theorem are satisfied we will show that $\del_0\neq0$.

For any $\ell$ set $F_\ell:=H^2(X;\ell)/\torsion$. If $b^+(X;\ell)=0$
let $D_\ell\subset F_\ell$ be
the subgroup generated by vectors of square $-1$.
Let $\hfl\subset F_\ell$ be the orthogonal complement of $D_\ell$, so that
\be{equation}\label{eqn:fell-dell-tfl}
F_\ell=D_\ell\oplus\hfl.
\end{equation}
By assumption there is a non-trivial $\ell$ such that $\hfl\neq0$ and 
$H^2(X;\ell)$ contains no element of order $4$. Fix such an $\ell$.
Note that there is a class $x\in\hfl$ with
$x^2\not\equiv0\mod4$. (Proof: Since $\hfl$ is unimodular we can find
elements $a,b\in\hfl$ with $a\cdot b=1$.
If $a^2,b^2\equiv0\mod4$ then $(a+b)^2\equiv2\mod4$.)
Let $k$ be the smallest (positive) integer $\not\equiv0\mod4$
such that there exists an $x\in\hfl$ with $x^2=-k$.

\subsection{Reduction to $\del(X)=0$}

By \Ref{eqn:del-indep} we have $\del(X)\le0$.
We will now reduce the remaining part of the proof to the case $\del(X)=0$.

\be{lemma}\label{lemma:taubpl}
Let $N$ be any compact, connected oriented smooth $4$--manifold, and
let $C$ be any embedded circle in $\itr(N)$ which represents a non-zero class
in $H_1(N;\z/2)$. Let $N'$ be obtained from $N$ by surgery on $C$. Then
\[\tau(N')+b^+(N')=\tau(N)+b^+(N).\]
\end{lemma}
Here $\tau$ is the invariant defined just before
Theorem~\ref{thm:main-theorem}.

\proof  Defining $\del(N)$ as in \Ref{eqn:deldef} we have
\be{align*}
b_1(N;\z/2)+\del(N)&=(b_1(N)+\tau(N))+(1-b_1(N)+b^+(N))\\
&=\tau(N)+b^+(N)+1.
\end{align*}
Now, $b_1(N;\z/2)$ drops by one by surgery on $C$, whereas
$\del(N)$ increases by one by surgery on {\em any} circle in $\itr(N)$.
(A highbrow proof of the latter statement applies the
excision principle for indices to the elliptic operator
$d^*+d^+:\Om^1\to\Om^0\oplus\Om^+$ on some 
close-up $V$ of $N$, recalling that the index of that operator is $-\del(V)$.)
\square

Every element of $H_1(X;\ell)$ can be represented by an embedded, oriented
circle $C$ in the interior of $X$ together with a trivialization of $\ell|_C$.
Set $\fd:=-\del(X)$ and let $X'$ be obtained from $X$ by performing
surgery on a collection of disjoint
oriented circles $C_1,\dots,C_{\fd}$ in $\itr(X)$
which, together with a trivialization of
$\ell$ over $\ti C:=\cup_jC_j$, represent a basis for $H_1(X;\ell)/\torsion$.
Let $\ell'\to X'$ be the bundle obtained by trivially extending
$\ell|_{X\setminus\ti C}$. Then $b_1(X';\ell')=0$, and there is a canonical
isomorphism $H^2(X;\ell)\to H^2(X';\ell')$ which induces an isomorphism
between the intersection forms. 
It follows from the long exact sequence
\[\cdots\to H_1(X;\ell)\oset{\cdot2}\to H_1(X;\ell)\to H_1(X;\z/2)\to\cdots\]
that the circles $C_j$ represent linearly independent classes in $H_1(X;\z/2)$,
so by Lemma~\ref{lemma:taubpl} the invariant $\tau+b^+$ takes the same value
on $X$ and $X'$.

We have shown that $X',\ell'$ satisfy all the hypotheses of the theorem, and
that $\ell',k$ satisfy the same minimality condition
as $\ell,k$. We may therefore
from now on assume that $b_1(X;\ell)=0=\del(X)$. This implies that
\be{equation}\label{eqn:ble3}
b:=b_1(X;\z/2)=\tau(X)+(b^+(X)+1)\le3,
\end{equation}
where we have used assumption \Ref{eqn:thmeqn} of the theorem.

\subsection{Choosing the sections}
\label{subsec:choosing-sec}

Let $W$ be the result of attaching a half-infinite
cylinder $[0,\infty)\times Y$ to $X$. We extend the bundle $\ell\to X$
to all of $W$ and, abusing notation, denote the new bundle also by $\ell$.
Choose a $c\in H^2(W;\ell)$ whose image in $F_\ell$ lies in $\hfl$ and such
that $c^2=-k$.
Define $\lla,K$ in terms of $\ell$
as in Section~\ref{sec:cohomology-local}
and let $L\to W$ be a Hermitian $K$--line bundle with $\tc(L)=c$. 
Then $E:=\lla\oplus L$ is an oriented, Euclidean $3$--plane bundle over $W$.

We will use the same notation for moduli spaces associated to $E$ as in
Section~\ref{sec:tw-red}.
Choose a Riemannian metric on $W$ which is on product form on the
end and which is generic 
as assumed in the beginning of Section~\ref{sec:tw-red}. 

Let $M^\lla_k$ be the set of all $[A]\in M_k$ such that $A$ preserves
a subbundle of $E$ isomorphic to $\lla$.
After perhaps perturbing the instanton equation as in
Section~\ref{sec:local-str-twr}
we may assume that every element of $M^\lla_k$ is a regular point in $M_k$.

We also add holonomy perturbations over the end of $W$ corresponding to
a small generic perturbation of the Chern--Simons functional over $Y$ (which
is in general needed to construct the Floer homology of $Y$), as well
as small holonomy perturbations obtained from a finite number of
thickened loops in $W$. (In order not to
obscure the main ideas, we usually ignore holonomy
perturbations in this paper.)

Let $M^\#_k$ be obtained from $M_k$
by deleting the interior of a small compact neighbourhood $N_\om$
of every $\om\in M^\lla_k$,
where $N_\om$ is as constructed in the proof
of Lemma~\ref{lemma:xi-link}. Let $M^-_k$ be the irreducible part of $M^\#_k$.

We are going to cut down $M^-_k$ to a $1$--manifold with boundary in the
following way.
For $i=1,\dots,2k-1$ let $Z'_i\subset W$ be a compact,
connected codimension~$0$ submanifold and $\ga_i:S^1\to Z'_i$ a loop.
Let $\Theta_{\ga_i}\to U_{\ga_i}$ be
the real line bundle associated to the double covering
$\Xi_{\ga_i}\to U_{\ga_i}$. Let $s'_i$ be a smooth section
of $\Theta_{\ga_i}$ over
the irreducible part of $U_{\ga_i}$ and set
\[\hmk:=\{\om\in M^-_k\st\text{$s'_i(\om|_{Z'_i})=0$ for $i=1,\dots,2k-1$}\}.\]
For generic sections $s'_i$ the space $\hmk$ will
be a smooth $1$--manifold with boundary (see \cite{Smale1}).
We will show that for a suitable choice of loops and
sections the manifold $\hat M_k$ will have an odd number of boundary
points and no ends coming from reducibles (i.e.\ points or circles in
$M\red$). We briefly outline how this will be achieved.

Consider the set
\[Q:=\{w\in H^1(W;\z/2)\st\text{$(\ga_i)^*w\neq0$ for
  $i=1,\dots,2k-1$}\}.\]
If $w_1(\lla)\in Q$ then, as we will see in
Section~\ref{subsec:ends-and-boundary-points}, $\hat M_k$ will have an
odd number of boundary points. If $w_1(\lla)$ is the {\em unique}
point in $Q$ then the sections $s'_i$ can be chosen so that $\hat M_k$
has no ends
associated to twisted reducibles in $M_k$. 
Note that $|Q|=1$ is indeed possible, since $b\le 3\le2k-1$.

To avoid ends of $\hat M_k$ associated to circles in $M\ared$ we choose
$Z'_{2j-1}=Z'_{2j}$ for $j=1,\dots,k-1$
and exploit the fact that the
direct sum of two real line bundles admits a non-vanishing section over
any circle (see Lemma~\ref{lemma:non-van2}); furthermore, we take the
$Z_j:=Z'_{2j-1}$, $j=1,\dots,k$ to be disjoint.

Finally, to arrange, in addition, that there are no ends in $\hat M_k$
coming from twisted reducibles in lower strata
$M_\ell$, $\ell<k$ we rotate the classes represented by the $\ga_i$ in a
suitable way.

We now make precise the choice of loops and sections.
Choose a basis $\{e_1,\dots,e_b\}$ for $H_1(W;\z/2)$ such that
$\la e_h,w_1(\lla)\ra=1$ for each $h$. Also choose 
\be{itemize}
\item disjoint compact, connected,
codimension~$0$ submanifolds $Z_1,\dots,Z_k$ of $W$,
\item two loops $\ga_{2j-1},\ga_{2j}$ in $Z_j$ for $j=1,\dots,k-1$,
\item a loop $\ga_{2k-1}$ in $Z_k$
\end{itemize}
such that $\ga_i$ represents $e_h$ when $i\equiv h\mod b$. For instance, $Z_k$
may be a closed tubular neighbourhood of an embedded circle in $W$, whereas
for $j=1,\dots,k-1$ one can take $Z_j$ to be an internal connected sum of two
such tubular neighourhoods.

We will write $U_i,\Xi_i,\cb_j$ instead of $U_{\ga_i},\Xi_{\ga_i},\cb(E|_{Z_j})$.
Let $\cb^*_j$ denote the irreducible part of $\cb_j$. Let $\Theta_i\to U_i$
be the real line bundle associated to the double covering $\Xi_i\to U_i$. 

For $j=1,\dots,k-1$ let $R_j\subset\cb_j$ be the image
of $M\red$ under the restriction map. We have observed that
$R_j$ is the disjoint union of finitely many circles and a finite set.
Note that, by Lemma~\ref{prop:red-in-Uga},
all these circles are contained in $U_{2j-1}\cap U_{2j}$.
Let $\hth_j$ be the direct sum of the
restrictions of $\Theta_{2j-1}$ and $\Theta_{2j}$ to $\cb^*_j$.
Let $s_j$ be a generic smooth section of $\hth_j$ which is
nowhere zero on $\cb^*_j\cap V_j$ for some neighbourhood
$V_j$ of the compact set $R_j\cap(U_{2j-1}\cup U_{2j})$
in $\cb_j$. The existence of 
sections of this kind follows from Lemma~\ref{lemma:non-van2}.
(The fact that $\cb_j$ is metrizable was pointed out in
Section~\ref{sec:conn-hol}.)

Let $R_k\subset\cb_k$ be the image of $M\tred$ under the restriction map.
Let $s_k$ be a generic smooth section of $\Theta_{2k-1}$ over $\cb^*_k$
which is nowhere zero on $\cb^*_k\cap V_k$ for some neighbourhood $V_k$
of $R_k\cap U_{2k-1}$ in $\cb_k$. The existence of 
sections of this kind follows from Lemma~\ref{lemma:non-van1}.

\subsection{Ends and boundary points}
\label{subsec:ends-and-boundary-points}

Set
\[\hat M_k:=\{\om\in M^-_k\st
\text{$s_j(\om|_{Z_j})=0$ for $j=1,\dots,k$}\}.\]
Modulo~$2$ the number of boundary points of the smooth $1$--manifold $\hmk$ is
\be{equation}\label{eqn:boundary-points}
\#\prtl\hmk\equiv\sum_{\om}
\la[\prtl N_\om],e(\Theta_\om)\ra\equiv|P_\ell|\equiv1\mod2,
\end{equation}
where $e$ denotes the Euler class with coefficients in $\z/2$ and
$\Theta_\om$ is the direct sum of the pull-backs of the line bundles
$\Theta_1,\dots,\Theta_{2k-1}$ to the boundary
$\prtl N_\om\approx\mathbb{RP}^{2k-1}$ of $N_\om$.

To prove the second congruence in \Ref{eqn:boundary-points}, note that
$\Theta_i$ pulls back to a non-trivial bundle
over each $\prtl N_\om$ by Lemma~\ref{lemma:xi-link}. Since the Euler class is
multiplicative under finite direct sums, we conclude that each term in the
sum in \Ref{eqn:boundary-points} is one.
The last congruence in \Ref{eqn:boundary-points} follows from
Proposition~\ref{prop:cardinality} because $|P_c|=1$ by the minimality
property of $k$, and $|2\ct_\ell|$ is odd since by assumption
$H^2(W;\ell)$ contains no element of order~$4$.

It remains to determine the ends of $\hmk$.
For any moduli space $M_{\al,d}$ with $\al$ irreducible set
\[\hat M_{\al,d}:=\{\om\in M_{\al,d}\st
\text{$s_j(\om|_{Z_j})=0$ for $j=1,\dots,k$}\}.\]

\be{prop}\label{prop:An-chain-conv}
Any sequence in $\hmk$ has a subsequence which either converges in
$\hmk$ or chain-converges to an element of
\[\hat M_{\al,2k-1}\times\check M(\al,\theta)\]
for some $\al\in\fl^*_Y$,
where $M(\al,\theta)$ is the one-dimensional moduli space over $\ry$
with limits $\al$ at $-\infty$ and $\theta$ at $\infty$, and
$\check M:=M/\R$.
\end{prop}

{\em Proof of proposition:} Let $\{[A_n]\}$ be a sequence in $\hmk$. After passing to
a subsequence we may assume that $\{[A_n]\}$
chain-converges weakly in the sense of \cite{D5}. Let $([A],x_1,\dots,x_q)$
be the weak limit over $W$, where $[A]\in M_{\al,d}$, $\al\in\fl_Y$,
and $x_1,\dots,x_q\in W$, $q\ge0$. We are going to show that $A$ must be
irreducible. First we establish the following lemma.

\be{lemma}\label{lemma:Aredcontr}
If $A$ is reducible then there is a $j\in\{1,\dots,k\}$ with
the following two properties:
\be{description}
\item[(i)]$Z_j$ contains none of the points $x_1,\dots,x_q$.
\item[(ii)]$[A|_{Z_j}]\in V_j$.
\end{description}
\end{lemma}

{\em Proof of lemma:} If $A$ is reducible, then $[A]\in M\red_{k-4s}$ for
some non-negative integer $s$. Observe that
\be{equation}\label{eqn:qsk}
q\le s<\frac k4<k-1.
\end{equation}
The second inequality holds because $k-4s\ge0$ by \Ref{eqn:fafa} and
we have chosen $k\not\equiv0\mod4$.
Hence there is certainly a $j<k$ satisfying (i).

{\bf Case 1:} $A$ Abelian. Then (ii) is satisfied for any $j<k$, so the
lemma holds in this case.

{\bf Case 2:} $A$ twisted reducible. Let $E=\lla'\oplus L'$ be the
splitting preserved by $A$, where $\lla'$ is a non-trivial real line bundle.

{\bf Case 2a:} $\lla'\approx\lla$. We will
show that this cannot occur.
Let $\ell'\subset\lla'$ be the lattice of vectors of integer length and set
$c':=\tc(L')\in H^2(W;\ell')$. Choose an isomorphism
$f:\ell'\oset\approx\to\ell$ and set $\zeta:=f_*c'\in H^2(W;\ell)$. Since
$[\zeta]_2=[c']_2=[c]_2$ (the last equality by
Proposition~\ref{prop:twisted-red-class}) it follows from the exact sequence
\Ref{eqn:coh-bockstein} that
there is an $x\in H^2(W;\ell)$ such that $\zeta=c+2x$.
For any $v\in H^2(W;\ell)$ let $\bar v$ be the image of $v$ in $F_\ell$ and
let $\hat v$ be the component of $\bar v$ in $\hfl$ with respect to the
splitting \Ref{eqn:fell-dell-tfl}. Since $\bar c\in\hfl$ by assumption, we have
$\hat\zeta=\bar c+2\hat x$, so $(\hat\zeta)^2\equiv c^2=-k\not\equiv0\mod4$.
Hence $-(\hat\zeta)^2\ge k$ by the minimality of $k$, so
\[k-4s=-\zeta^2\ge-(\hat\zeta)^2\ge k.\]
Thus, $s=0$ and $[A]\in M_k$. It follows that the sequence
$\{[A_n]\}$ converges in $M_k$ (see \cite{D5}). Since $M^\#_k$ is a closed
subset of $M_k$, we must have $[A]\in M^\#_k$. This is a contradiction, since
$M^\#_k$ was obtained from $M_k$ by deleting neighbourhoods of all twisted
reducibles preserving a line bundle isomorphic to $\lla$.

{\bf Case 2b:} $\lla'\not\approx\lla$. Then $b\ge2$. 

{\bf Case 2b1:} $b=2$. Then for $h=1$ or $2$ we have
\[1=\la e_h,w_1(\lla)+w_1(\lla')\ra=1+\la e_h,w_1(\lla')\ra,\]
so $\la e_h,w_1(\lla')\ra=0$. As observed in the beginning of the
proof we can find a $j<k$ satisfying (i).
For $i=2j-1$ or $2j$ the loop
$\ga_i$ represents $e_h$, in which case $(\ga_i)^*\lla'$ is
trivial. This in turn implies $[A|_{Z_j}]\in U_i$ by
Proposition~\ref{prop:red-in-Uga}, so $[A|_{Z_j}]\in V_j$.

{\bf Case 2b2:} $b=3$. Set
\[m:=\left[\frac{k-1}3\right].\]

{\bf Case 2b2a:} $m=0$. Then
$k\le3$, so $q=0$ by \Ref{eqn:qsk}. The same argument as in the case
$b=2$ shows that $(\ga_i)^*\lla'$ is trivial for some $i\in\{1,2,3\}$,
hence $[A|_{Z_j}]\in V_j$ for $j=1$ or $2$.

{\bf Case 2b2b:} $m\ge1$. Choose $h$ with $\la e_h,w_1(\lla')\ra=0$.
Then for at least $2m$ integers
$j\in\{1,\dots,k-1\}$ one of the loops $\ga_{2j-1}$ or $\ga_{2j}$
will represent $e_h$, in which case (ii) holds. Because
\[q<\frac k4<2m,\]
one can choose $j$ such that (i) holds as well.\square \\

\be{lemma} $A$ is irreducible.
\end{lemma}

{\em Proof of lemma:} Assume to the contrary that $A$ were reducible,
and let $j$ satisfy the two properties of Lemma~\ref{lemma:Aredcontr}.
Then
\[[A_n|_{Z_j}]\to[A|_{Z_j}]\quad\text{in $\cb_j$ as $n\to\infty$}.\]
But $V_j$ is open, so for sufficently large $n$
we have $[A_n|_{Z_j}]\in\cb^*_j\cap V_j$ and therefore $s_j(A_n|_{Z_j})\neq0$.
This contradicts $[A_n]\in\hmk$ and the lemma is proved.\square

We can now complete the proof of the proposition. First suppose $[A]\in
M_k$, which implies $q=0$. Then
$\{[A_n]\}$ converges in $M_k$, so
$[A]\in M^\#_k$. But $[A]$ is irreducible, so $[A]\in\hmk$.

Now suppose $[A]\not\in M_k$. Then
\be{equation}\label{eqn:dlemin}
d\le\min(2k-1,2k-8q).
\end{equation}
Set
\[J:=\{j\in\{1,\dots,k-1\}\st
\text{$Z_j$ contains none of the points $x_1,\dots,x_q$}\}.\]
Then $s_j(A|_{Z_j})=0$ for every $j\in J$. Since the sections $s_j$
are generic, we must have $2|J|\le d$, where $|J|$ denotes the
cardinality of the set $J$. Combining this with \Ref{eqn:dlemin}
yields
\[2|J|\le2k-8q.\]
Setting $t:=k-1-|J|$ we deduce
\[4q\le t+1\le q+1,\]
so $q=0$. Hence $s_j(A|_{Z_j})=0$ for $j=1,\dots,k$, so
$d\ge2k-1$. Combining this with \Ref{eqn:dlemin} we obtain
$d=2k-1$. This is only possible when $\al$ is irreducible, so the
proposition is proved.\square

We can now complete the proof of the theorem.
An argument similar to the proof of
Proposition~\ref{prop:An-chain-conv}
shows that $\hat M_{\al,2k-1}$ is compact, hence a finite
set (since it is $0$--dimensional). By gluing theory the number of ends of
$\hmk$ is $\del h$, where
\[h:=\sum_\al[\#\hat M_{\al,2k-1}]\al\in\CF^4(Y)\otimes\z/2.\]
The proof of the proposition applied
to moduli spaces $M_{\beta,2k}$ with $\beta$ irreducible shows that $h$ is
a cocycle.

Since the $1$--manifold $\hmk$ has an odd number of boundary points, it must
also have an odd number of ends, so $\del_0([h])=1$.
\square

\noindent\textsc{Centre for Quantum Geometry of Moduli Spaces,\\
Ny Munkegade 118,\\
DK-8000 Aarhus C,\\
Denmark}
\\ \\
Email:\ froyshov@imf.au.dk


\begin{thebibliography}{10}

\bibitem{Bredon}
G.~E. Bredon.
\newblock {\em Sheaf Theory. Second edition}.
\newblock Springer, 1997.

\bibitem{Dold-Whitney}
A.~Dold and H.~Whitney.
\newblock Classification of oriented sphere bundles over a 4--complex.
\newblock {\em Annals of Math.}, 69(3):667--677, 1959.

\bibitem{D1}
S.~K. Donaldson.
\newblock An application of gauge theory to four dimensional topology.
\newblock {\em J.~Diff. Geom.}, 18:279--315, 1983.

\bibitem{D2}
S.~K. Donaldson.
\newblock The orientation of {Yang--Mills} moduli spaces and $4$--manifold
  topology.
\newblock {\em J.~Diff. Geom.}, 26:397--428, 1987.

\bibitem{D5}
S.~K. Donaldson.
\newblock {\em Floer Homology Groups in {Yang--Mills} Theory}.
\newblock Cambridge University Press, 2002.

\bibitem{DK}
S.~K. Donaldson and P.~B. Kronheimer.
\newblock {\em The Geometry of Four-Manifolds}.
\newblock Oxford University Press, 1990.

\bibitem{FS2}
R.~Fintushel and R.~J. Stern.
\newblock Definite 4--manifolds.
\newblock {\em J.~Diff. Geom.}, 28:133--141, 1988.

\bibitem{F1}
A.~Floer.
\newblock An instanton invariant for 3--manifolds.
\newblock {\em Comm. Math. Phys.}, 118:215--240, 1988.

\bibitem{FU}
D.~S. Freed and K.~K. Uhlenbeck.
\newblock {\em Instantons and Four-Manifolds}.
\newblock MSRI Publications. Springer-Verlag, second edition, 1991.

\bibitem{Freedman-Quinn}
M.~Freedman and F.~Quinn.
\newblock {\em Topology of $4$--manifolds}, volume~39 of {\em Princeton
  Mathematical Series}.
\newblock Princeton University Press, 1990.

\bibitem{FHMT}
S.~Friedl, I.~Hambleton, P.~Melvin, and P.~Teichner.
\newblock Non-smoothable four-manifolds with infinite cyclic fundamental group.
\newblock {\em Int. Math. Res. Not.}, (11), 2007.

\bibitem{Fr3}
K.~A. Fr{\o}yshov.
\newblock Equivariant aspects of {Yang--Mills} {Floer} theory.
\newblock {\em Topology}, 41:525--552, 2002.

\bibitem{Fr13}
K.~A. Fr{\o}yshov.
\newblock {\em {Compactness and gluing theory for monopoles}}, volume~15 of
  {\em {Geometry \&\ Topology Monographs}}.
\newblock {Geometry \&\ Topology Publications}, 2008.

\bibitem{Fur2}
M.~Furuta.
\newblock Monopole equation and the $\frac{11}8$--conjecture.
\newblock {\em Math. Res. Lett.}, 8:279--291, 2001.

\bibitem{Hambleton1}
I.~Hambleton.
\newblock Intersection forms, fundamental groups and 4-manifolds.
\newblock In {\em {Proceedings of G\"okova Geometry-Topology Conference 2008}},
  pages 137--150, 2009.

\bibitem{Hatcher1}
A.~Hatcher.
\newblock {\em Algebraic Topology}.
\newblock Cambridge University Press, 2002.

\bibitem{Kelley}
J.~L. Kelley.
\newblock {\em General Topology}.
\newblock Van Nostrand, 1955.

\bibitem{KM3}
P.~B. Kronheimer and T.~S. Mrowka.
\newblock Embedded surfaces and the structure of {Donaldson's} polynomial
  invariants.
\newblock {\em J.~Diff. Geom.}, 41:573--734, 1995.

\bibitem{Lang1}
S.~Lang.
\newblock {\em Differential Manifolds. Second edition}.
\newblock Springer-Verlag, 1985.

\bibitem{nakamura1}
N.~Nakamura.
\newblock {$\text{Pin}^-(2)$--monopole equations and intersection forms with
  local coefficients of $4$--manifolds}.
\newblock arXiv:1009.3624.

\bibitem{Smale1}
S.~Smale.
\newblock {An infinite dimensional version of Sard's theorem}.
\newblock {\em Am. J. Math.}, 87:861--866, 1965.

\bibitem{Spanier}
E.~H. Spanier.
\newblock {\em Algebraic Topology. Corrected reprint}.
\newblock Springer, 1981.

\bibitem{Spanier2}
E.~H. Spanier.
\newblock Singular homology and cohomology with local coefficients and duality
  for manifolds.
\newblock {\em Pacific J. Math.}, 160:165--200, 1993.

\bibitem{Teleman1}
A.~Teleman.
\newblock Harmonic sections in sphere bundles, normal neighborhoods of
  reduction loci, and instanton moduli spaces on definite $4$--manifolds.
\newblock {\em Geometry \&\ Topology}, 11:1681--1730, 2007.

\end{thebibliography}
\end{document}